\documentclass[leqno,12pt]{amsart}
\usepackage{amsfonts}
\usepackage{amsmath,amssymb,amsthm}

\newtheorem{thm}{Theorem}[section]

\newtheorem{prop}[thm]{Proposition}
\theoremstyle{definition}

\theoremstyle{remark}

\begin{document}

\title[Meridian Surfaces of Parabolic Type in the
Minkowski 4-space] {Meridian Surfaces of Parabolic Type in the
Four-dimensional Minkowski Space}

\author{Georgi Ganchev and Velichka Milousheva}
\address{Institute of Mathematics and Informatics, Bulgarian Academy of Sciences,
Acad. G. Bonchev Str. bl. 8, 1113 Sofia, Bulgaria}
\email{ganchev@math.bas.bg}
\address{Institute of Mathematics and Informatics, Bulgarian Academy of Sciences,
Acad. G. Bonchev Str. bl. 8, 1113, Sofia, Bulgaria;   "L. Karavelov"
Civil Engineering Higher School, 175 Suhodolska Str., 1373 Sofia,
Bulgaria} \email{vmil@math.bas.bg}


\begin{abstract}
We construct a special class of spacelike surfaces in the
Minkowski 4-space which are one-parameter systems of meridians of
the rotational hypersurface with lightlike axis and call these
surfaces meridian surfaces of parabolic type. They are analogous
to the meridian surfaces of elliptic or hyperbolic type. Using the
invariants of these surfaces we give the complete classification
of the meridian surfaces of parabolic type with constant Gauss
curvature or constant mean curvature. We also classify the Chen
meridian surfaces  of parabolic type and the meridian surfaces  of
parabolic type with parallel normal bundle.
\end{abstract}

\maketitle

\section{Introduction}

A fundamental problem of the contemporary
differential geometry of surfaces in  the Euclidean space $\mathbb{R}^n$ or the
pseudo-Euclidean  space $\mathbb{R}^n_k$ is the investigation of the basic
invariants characterizing the surfaces. Our aim is to investigate various important classes of surfaces in the
four-dimensional Minkowski space $\mathbb{R}^4_1$ characterized by
conditions on their invariants.

In \cite {GM5} we developed a local theory of spacelike surfaces in  $\mathbb{R}^4_1$ based on the introducing of an invariant linear map $\gamma$ of Weingarten-type  in the tangent plane at any point of the surface.
The map $\gamma$  generates two invariant functions $k = \det \gamma$ and $\varkappa= -\displaystyle{ \frac{1}{2}}\, \mathrm{tr} \gamma$. It turns out that the invariant  $\varkappa$ is the curvature of the normal connection of the surface.
The existence of principal lines at each point of a spacelike surface in $\mathbb{R}^4_1$ allows us to introduce a
geometrically determined moving frame field at each point of the
surface. Writing derivative formulas for this frame
field, we obtained eight invariant functions  $\gamma_1, \,
\gamma_2, \, \nu_1,\, \nu_2, \, \lambda, \, \mu, \, \beta_1,
\beta_2$ and proved a fundamental theorem of Bonnet-type, stating
that these eight invariants under some natural conditions
determine the surface up to a rigid motion in $\mathbb{R}^4_1$.

The basic geometric classes of surfaces in $\mathbb{R}^4_1$   are
characterized by conditions on these invariant functions. For
example, Chen surfaces are characterized by the condition $\lambda = 0$, minimal surfaces are determined by the equality
$\nu_1 + \nu_2 = 0$, surfaces with flat normal connection are  described by $\nu_1 = \nu_2$,  and surfaces with parallel normal bundle are characterized by $\beta_1 = \beta_2 = 0$.

In the four-dimensional Minkowski space $\mathbb{R}^4_1$  there are three types of
rotational hypersurfaces  -  rotational hypersurfaces with
timelike axis, with spacelike axis, and with lightlike axis. In
\cite{GM6} we  constructed special families of two-dimensional
spacelike surfaces lying
 on rotational hypersurfaces in $\mathbb{R}^4_1$  with
timelike or spacelike axis and called them meridian surfaces of
elliptic or hyperbolic type, respectively. We found all
marginally trapped meridian surfaces of elliptic or hyperbolic
type.
 In \cite{GM-Ell-Hyp} we found  the geometric invariant
functions $\gamma_1, \, \gamma_2, \, \nu_1,\, \nu_2, \, \lambda,
\, \mu, \, \beta_1, \beta_2$ of the meridian surfaces of elliptic
or hyperbolic type and classified those of them with constant
Gauss curvature or constant mean curvature. We also gave the
complete classification of the Chen meridian surfaces  of elliptic
or hyperbolic type  and the meridian surfaces of elliptic or
hyperbolic type with parallel normal bundle.

In \cite{GM7} we used the idea from the elliptic and
hyperbolic case to construct families of two-dimensional
spacelike surfaces lying
 on a rotational hypersurface in $\mathbb{R}^4_1$  with
lightlike axis.  We called these surfaces \emph{meridian surfaces of
parabolic type}.  We found all
marginally trapped meridian surfaces of parabolic type.

In the present paper we study meridian surfaces of parabolic type in $\mathbb{R}^4_1$ and find the invariant functions $\gamma_1, \,
\gamma_2, \, \nu_1,\, \nu_2, \, \lambda, \, \mu, \, \beta_1,
\beta_2$ of these surfaces. Using the invariants we  classify completely  the meridian surfaces  of parabolic type with
constant Gauss curvature (Theorem \ref{T:Gauss curvature}), with constant mean curvature (Theorem \ref{T:mean curvature}), and  with constant
invariant $k$ (Theorem \ref{T:constant k}). In Theorem \ref{T:Chen}  we classify  the Chen meridian surfaces of parabolic type
and in Theorem \ref{T:parallel}  we give the classification of the meridian surfaces of parabolic type
with parallel normal bundle.

\section{Invariants of meridian surfaces of parabolic type}

We consider the  four-dimensional Minkowski space  $\mathbb{R}^4_1$  endowed with the metric
$\langle , \rangle$ of signature $(3,1)$. 
A surface $M^2: z = z(u,v), \, \, (u,v) \in {\mathcal D}$
(${\mathcal D} \subset \mathbb{R}^2$) in $\mathbb{R}^4_1$ is said to be
\emph{spacelike} if $\langle , \rangle$ induces  a Riemannian
metric $g$ on $M^2$.
Denote by $\nabla'$ and $\nabla$ the Levi Civita connections on $\mathbb{R}^4_1$ and $M^2$, respectively.
Let $x$ and $y$ be vector fields tangent to $M^2$ and $\xi$ be a normal vector field.
The formulas of Gauss and Weingarten give the decompositions of the vector fields $\nabla'_xy$ and
$\nabla'_x \xi$ into tangent and normal components \cite{Chen1}:
$$\begin{array}{l}
\vspace{2mm}
\nabla'_xy = \nabla_xy + \sigma(x,y);\\
\vspace{2mm}
\nabla'_x \xi = - A_{\xi} x + D_x \xi,
\end{array}$$
which define the second fundamental tensor $\sigma$, the normal connection $D$
and the shape operator $A_{\xi}$ with respect to $\xi$.
The mean curvature vector  field $H$ of $M^2$ is defined as
$H = \displaystyle{\frac{1}{2}\,  \mathrm{tr}\, \sigma}$.

Studying spacelike surfaces in $\mathbb{R}^4_1$ whose mean
curvature vector at any point is a non-zero spacelike vector or
timelike vector, on the base of the principal lines we introduced
a geometrically determined orthonormal frame field $\{x,y,b,l\}$ at each point of such a surface  \cite{GM5}.
The tangent vector fields $x$ and $y$ are collinear with  the principal directions,
the normal vector field $b$ is collinear with the mean curvature vector field $H$.
Writing derivative formulas
of Frenet-type for this frame field, we obtained eight invariant
functions $\gamma_1, \, \gamma_2, \, \nu_1,\, \nu_2, \, \lambda, \, \mu,
\, \beta_1, \beta_2$, which determine the surface up to a rigid motion in $\mathbb{R}^4_1$.
These invariants are determined by the geometric frame field $\{x,y,b,l\}$ as follows
\begin{equation}
\begin{array}{l} \label{E:Eq1}
\vspace{2mm}
\nu_1 = \langle \nabla'_xx, b\rangle, \qquad \nu_2 = \langle \nabla'_yy, b\rangle, \qquad \lambda = \langle \nabla'_xy, b\rangle,
\qquad \mu = \langle \nabla'_xy, l\rangle,\\
\vspace{2mm}
\gamma_1 =  \langle \nabla'_xx, y\rangle, \qquad  \gamma_2 = \langle \nabla'_yy, x\rangle, \qquad \beta_1 = \langle \nabla'_xb, l\rangle, \qquad
\beta_2 = \langle \nabla'_yb, l\rangle.
\end{array}
\end{equation}

The invariants $k$, $\varkappa$, and the Gauss curvature $K$ of
$M^2$ are expressed by the functions $\nu_1, \nu_2, \lambda, \mu$
as follows:
\begin{equation} \notag
k = - 4\nu_1\,\nu_2\,\mu^2, \quad \quad \varkappa = (\nu_1-\nu_2)\mu, \quad \quad
K = \varepsilon (\nu_1\,\nu_2- \lambda^2 + \mu^2),
\end{equation}
where $\varepsilon = sign \langle H, H \rangle$.

\vskip 2mm In the present section we give the construction of meridian surfaces of parabolic type and find their invariant functions
$\gamma_1, \, \gamma_2, \, \nu_1,\, \nu_2, \, \lambda, \, \mu, \, \beta_1, \beta_2$.

Let $\{e_1, e_2, e_3, e_4\}$ be the standard orthonormal frame in the Minkowski space
$\mathbb{R}^4_1$, i.e.  $\langle e_1, e_1 \rangle  =
\langle e_2, e_2 \rangle  = \langle e_3, e_3 \rangle  = 1, \, \langle e_4, e_4 \rangle  = -1$. We
denote  $\displaystyle{\xi_1= \frac{e_3 + e_4}{\sqrt{2}}},\,\, \displaystyle{\xi_2= \frac{ - e_3 + e_4}{\sqrt{2}}}$ and consider the
 pseudo-orthonormal base $\{e_1, e_2,
\xi_1, \xi_2 \}$  of $\mathbb{R}^4_1$.
Note that $\langle\xi_1, \xi_1 \rangle =0$, $\langle \xi_2, \xi_2
\rangle =0$, $\langle \xi_1, \xi_2 \rangle = -1$. The rotational
hypersurface with lightlike axis  in $\mathbb{R}^4_1$ can be parameterized by
$$\mathcal{M}''': Z(u,w^1,w^2) =  f(u)\, w^1 \cos w^2 \,e_1 +  f(u)\, w^1 \sin w^2 \,e_2+ \left(f(u) \frac{(w^1)^2}{2} + g(u)\right) \xi_1 + f(u) \,\xi_2,$$
where $f = f(u), \,\, g = g(u)$ are smooth functions, defined in
an interval $I \subset \mathbb{R}$, such that $- f'(u)g'(u) >0$, $f(u)>0,
\,\, u \in I$.

Let $w^1 = w^1(v), \, w^2=w^2(v), \,\, v \in J, \,J \subset \mathbb{R}$
and assume that $(\dot{w}^1)^2 + (\dot{w}^2)^2 \neq 0$.
 We consider the surface $\mathcal{M}'''_m$ in $\mathbb{R}^4_1$ defined as follows
\begin{equation} \label{E:Eq-1}
\mathcal{M}'''_m: z(u,v) = Z(u,w^1(v),w^2(v)),
\end{equation}
where $u \in I, \, v \in J.$ The surface $\mathcal{M}'''_m$,
defined by \eqref{E:Eq-1}, is a one-parameter system of meridians
of the rotational hypersurface $\mathcal{M}'''$  with lightlike
axis. We  call $\mathcal{M}'''_m$ a \textit{meridian surface
of parabolic type}.

\vskip 2mm

Without loss of generality we assume that $w^1 = \varphi(v), \,
w^2=v$. Then the surface $\mathcal{M}'''_m$ is parameterized as
follows:
\begin{equation} \label{E:Eq-2}
\mathcal{M}'''_m: z(u,v) = f(u)\, \varphi(v) \cos v \,e_1 + f(u)\,
\varphi(v) \sin v \,e_2+ \left(f(u) \frac{(\varphi(v))^2}{2} +
g(u)\right)\xi_1 + f(u) \,\xi_2.
\end{equation}

The parametric $u$-lines of   $\mathcal{M}'''_m$ are curves congruent in $\mathbb{R}^4_1$ and
the curvature of each $u$-line is $\displaystyle{ \frac{f' g'' - g' f''}{(-2 f' g')^{\frac{3}{2}}}}$. These curves are the meridians of  $\mathcal{M}'''_m$.
We denote $\kappa_m(u) = \displaystyle{ \frac{f' g'' - g' f''}{(-2 f' g')^{\frac{3}{2}}}}$.

For each  $u = u_0 = const$  the curvature of the corresponding parametric  $v$-line is
$\varkappa_{c_v} = \displaystyle{\frac{\varphi \ddot{\varphi} - 2 \dot{\varphi}^2 - \varphi^2 }{ a(\dot{\varphi} ^2 + \varphi^2)^{\frac{3}{2}}}}$, where $a =f(u_0)$.
Let us denote $\kappa(v) = \displaystyle{\frac{\varphi \ddot{\varphi} - 2 \dot{\varphi}^2 - \varphi^2 }{(\dot{\varphi} ^2 + \varphi^2)^{\frac{3}{2}}}}$.
Then, the curvature of the  $v$-line $u = u_0$  is expressed as
$\varkappa_{c_v} = \displaystyle{\frac{1}{a}\, \kappa(v)}$ \cite{GM7}.

The tangent vector fields of $\mathcal{M}'''_m$ are:
\begin{equation} \label{E:Eq-6}
\begin{array}{l}
\vspace{2mm}
z_u = \displaystyle{f' \varphi \cos v \,e_1 + f'
\varphi \sin v \,e_2+ \left(f' \frac{\varphi^2}{2} +
g'\right)\xi_1 + f' \,\xi_2};\\
\vspace{2mm}
z_v = f( \dot{\varphi} \cos v - \varphi \sin v)\,e_1 + f (
\dot{\varphi} \sin v + \varphi \cos v) \,e_2+ f \varphi \dot{\varphi} \, \xi_1.
\end{array}
\end{equation}
So, the coefficients of the first fundamental form of $\mathcal{M}'''_m$
are
$$E = - 2 f'(u) g'(u); \quad F = 0; \quad G = f^2(u) (\dot{\varphi}^2(v) + \varphi^2(v)).$$
Note that the first fundamental form is positive definite, since $- f' g' >0$.

Without loss of generality we assume that $ - 2 f'(u) g'(u) = 1$, i.e. the meridians are parameterized by the arc-length.
Then $\kappa_m(u) = \displaystyle{\frac{f''(u)}{f'(u)}}$.

Denote
$X = z_u,\,\,Y = \displaystyle{\frac{z_v}{f \sqrt{\dot{\varphi}^2 + \varphi^2}}}$ and
consider the following orthonormal normal frame field:
\begin{equation} \label{E:Eq-7}
\begin{array}{l}
\vspace{2mm}
n_1 = \displaystyle{\frac{1}{\sqrt{\dot{\varphi}^2 + \varphi^2}}  \left((\dot{\varphi} \sin v + \varphi \cos v)\,e_1 + (-\dot{\varphi} \cos v + \varphi \sin v)\,e_2+
 \varphi^2  \, \xi_1\right)};\\
\vspace{2mm}
n_2 = \displaystyle{ f' \left( \varphi \cos v \,e_1 + \varphi \sin v \,e_2 + \frac{f' \varphi^2 - 2g'}{f'} \, \xi_1 +
 \xi_2\right)}.
\end{array}
\end{equation}

Thus we obtain a frame field $\{X, Y, n_1, n_2\}$ of $\mathcal{M}'''_m$, such that $\langle n_1, n_1 \rangle =1$,
$\langle n_2, n_2 \rangle =- 1$, $\langle n_1, n_2 \rangle =0$.

Calculating the second partial derivatives of $z(u,v)$ and taking into account \eqref{E:Eq-7} we get the following derivative formulas:

\begin{equation} \label{E:Eq-5}
\begin{array}{ll}
\vspace{2mm} \nabla'_XX = \qquad \qquad \qquad \qquad
- \kappa_m\,n_2; & \qquad
\nabla'_X n_1 = 0;\\
\vspace{2mm} \nabla'_XY = 0;  & \qquad
\nabla'_Y n_1 = \displaystyle{\quad \quad \, - \frac{\kappa}{f}\,Y};\\
\vspace{2mm} \nabla'_YX = \quad\quad
\quad\displaystyle{\frac{f'}{f}}\,Y;  & \qquad
\nabla'_X n_2 = \kappa_m \,X;\\
\vspace{2mm} \nabla'_YY = \displaystyle{- \frac{f'}{f}\,X \quad\quad +
\frac{\kappa}{f}\,n_1 - \frac{f'}{f} \, n_2}; & \qquad \nabla'_Y
n_2 = \displaystyle{ \quad \quad \;\;  \frac{f'}{f}\,Y},
\end{array}
\end{equation}

The invariants $k$ and $\varkappa$ of the meridian surface of parabolic type are given by the following formulas:
$$k = - \frac{\kappa_m^2(u) \, \kappa^2(v)}{f^2(u)}; \qquad \varkappa = 0.$$

Since $\varkappa$ is the curvature of the normal connection, from the equality $\varkappa = 0$ we get the following result.

\begin{prop}
The meridian surface of parabolic type $\mathcal{M}'''_m$, defined by \eqref{E:Eq-2}, is a surface with
flat normal connection.
\end{prop}

Taking into account
\eqref{E:Eq-5} and using that $\kappa_m = \displaystyle{\frac{f''}{f'}}$, we find the Gauss
curvature $K$ and the mean curvature vector field $H$ of $\mathcal{M}'''_m$:
\begin{equation} \label{E:Eq-11}
K = \displaystyle{- \frac{f''(u)}{f(u)}};\\
\end{equation}

\begin{equation} \label{E:Eq-12}
H = \displaystyle{ \frac{\kappa(v)}{2f(u)}\,\, n_1 - \frac{f(u) f''(u) +
f'^2(u)}{2 f(u) f'(u)}  \, \, n_2}.
\end{equation}

We distinguish the following  three cases (see \cite{GM7}):

\vskip 2mm I.
$\kappa(v) = 0$. In this case $n_1 = const$ and $\mathcal{M}'''_m$
lies in the hyperplane $\mathbb{R}^3_1$ of $\mathbb{R}^4_1$ orthogonal to $n_1$,
i.e. $\mathcal{M}'''_m$ lies in $\mathbb{R}^3_1 = \mathrm{span} \{x,y,n_2\}$.

\vskip 2mm
II. $\kappa_m(u) = 0$. In this case  $\mathcal{M}'''_m$ is a developable ruled surface in $\mathbb{R}^4_1$.

\vskip 2mm
III. $\kappa_m(u) \, \kappa(v) \neq 0$.

\vskip 2mm
In the first two cases the surface $\mathcal{M}'''_m$ consists of flat points, i.e  $k = \varkappa = 0$. It is known that surfaces consisting of flat points either lie in a hyperplane of $\mathbb{R}^4_1$ or are developable ruled surfaces. So, we consider the  third (general) case,
 i.e. we assume that $\kappa_m \neq 0$ and $\kappa \neq 0$.

The mean curvature vector field $H$ of $\mathcal{M}'''_m$ is expressed by formula \eqref{E:Eq-12}.
Since $\kappa \neq 0$ the surface $\mathcal{M}'''_m$ is non-minimal, i.e. $H \neq 0$. Recall that a spacelike surface in $\mathbb{R}^4_1$ is called marginally trapped if $H \neq 0$ and $\langle H, H \rangle = 0$.
The marginally trapped meridian surfaces of parabolic type are described in \cite{GM7}.
So, here we consider the case  $\langle H, H \rangle \neq 0$.

The orthonormal frame field $\{X,Y, n_1, n_2\}$, defined above is not the  geometric frame field  of the surface, since $X$ and $Y$ are not principal tangents.
The principal tangents of $\mathcal{M}'''_m$  are determined by
\begin{equation}\notag
x = \displaystyle{\frac{X+Y}{\sqrt{2}}}; \qquad   y = \displaystyle{\frac{- X + Y}{\sqrt{2}}}.
\end{equation}
In the case $\langle H, H \rangle > 0$, i.e. $\kappa^2 f'^2 - (f
f'' + f'^2)^2 > 0$, the geometric normal frame field  $\{b,l\}$ is
given by
\begin{equation}\notag
\begin{array}{l}
\vspace{2mm}
b = \displaystyle{\frac{1}{\sqrt{\kappa^2 f'^2 - (f f'' + f'^2)^2}} \left( \kappa f'\, n_1 - (f f'' + f'^2)\,n_2 \right)}; \\
\vspace{2mm} l = \displaystyle{\frac{1}{\sqrt{\kappa^2 f'^2 - (f f'' +
f'^2)^2}} \left( (f f'' + f'^2) \, n_1 -  \kappa f'\,n_2 \right)}.
\end{array}
\end{equation}
In this case the normal vector fields $b$ and $l$ satisfy $\langle b, b \rangle = 1$, $\langle b, l \rangle = 0$, $\langle l, l \rangle = -1$.

In the case $\langle H, H \rangle < 0$,  i.e. $\kappa^2 f'^2 - (f f'' + f'^2)^2 < 0$,  the geometric normal frame field  $\{b,l\}$ is given by
\begin{equation}\notag
\begin{array}{l}
\vspace{2mm}
b = \displaystyle{- \frac{1}{\sqrt{(f f'' + f'^2)^2 - \kappa^2 f'^2 }} \left(\kappa f'\, n_1 - (f f'' + f'^2)\,n_2 \right)}; \\
\vspace{2mm} l = \displaystyle{\frac{1}{\sqrt{(f f'' + f'^2)^2 - \kappa^2
f'^2}} \left( -(f f'' + f'^2) \, n_1 +  \kappa f'\,n_2 \right)}.
\end{array}
\end{equation}
In this case we have $\langle b, b \rangle = - 1$, $\langle b, l \rangle = 0$, $\langle l, l \rangle = 1$.

Using the  geometric frame field $\{x,y, b, l\}$ of $\mathcal{M}'''_m$  and derivative formulas
\eqref{E:Eq-5},
we obtain that the geometric  invariant functions   of  $\mathcal{M}'''_m$ are expressed by the formulas:

\begin{equation} \label{E:Eq7}
\begin{array}{l}
\vspace{2mm}
\gamma_1 = - \gamma_2 = \displaystyle{ \frac{f'}{\sqrt{2}f}};\\
\vspace{2mm}
\nu_1  = \nu_2 = \displaystyle{ \frac{\sqrt{\varepsilon (\kappa^2 f'^2 - (f f'' + f'^2)^2)}}{2f f'} };\\
\vspace{2mm}
\lambda = \displaystyle{ \varepsilon \frac{\kappa^2 f'^2 + f^2 f''^2 - f'^4}{2f f' \sqrt{\varepsilon (\kappa^2 f'^2 - (f f'' + f'^2)^2)}}}; \\
\vspace{2mm}
\mu = \displaystyle{\frac{\kappa f''}{\sqrt{\varepsilon (\kappa^2 f'^2 - (f f'' + f'^2)^2)}}}; \\
\vspace{2mm}
\beta_1 = \displaystyle{\frac{ - f'^2}{\sqrt{2} \varepsilon(\kappa^2 f'^2 - (f f'' + f'^2)^2)}} \left( \kappa\, \frac{d}{du}\left( \frac{f f'' + f'^2}{f'} \right) -  \frac{d}{dv}(\kappa) \, \frac{f f'' + f'^2}{f f' \sqrt{\dot{\varphi} ^2 + \varphi^2}} \right); \\
\vspace{2mm}
\beta_2 = \displaystyle{\frac{f'^2}{\sqrt{2} \varepsilon(\kappa^2 f'^2 - (f f'' + f'^2)^2)}} \left( \kappa\, \frac{d}{du}\left( \frac{f f'' + f'^2}{f'} \right) +  \frac{d}{dv}(\kappa) \, \frac{f f'' + f'^2}{f f' \sqrt{\dot{\varphi} ^2 + \varphi^2}} \right),
\end{array}
\end{equation}
where $\varepsilon = sign \langle H, H \rangle$.

\vskip 2mm
In the following sections, using the invariants of the meridian surface $\mathcal{M}'''_m$,
we shall describe and classify some special classes of meridian surfaces of parabolic type.

\section{Meridian surfaces  of parabolic type with constant Gauss curvature}

The study of surfaces with constant Gauss curvature is one of the
main topics in differential geometry. Surfaces
with constant Gauss curvature in Minkowski space have drawn the interest of many
geometers, see for example \cite{GalMarMil},  \cite{Lop}, and the references therein.

The Gauss curvature of a meridian surface  of parabolic type $\mathcal{M}'''_m$ depends only on the meridian curve $m$ and is expressed by  formula
\eqref{E:Eq-11}. The following theorem describes the  meridian surfaces of parabolic type with constant non-zero Gauss curvature.

\begin{thm} \label{T:Gauss curvature}
Let $\mathcal{M}'''_m$ be a meridian surface of parabolic  type from the general class.
Then $\mathcal{M}'''_m$ has constant non-zero Gauss curvature $K$ if and only
if the meridian $m$ is given by
\begin{equation} \label{E:Eq-Th}
\begin{array}{ll}
\vspace{2mm}
f(u) = \alpha \cos \sqrt{K} u + \beta \sin \sqrt{K} u, & \textrm{if} \quad K >0;\\
\vspace{2mm} f(u) = \alpha \cosh \sqrt{-K} u + \beta \sinh
\sqrt{-K} u, & \textrm{if} \quad K <0,
\end{array}
\end{equation}
where $\alpha$ and $\beta$ are constants, $g(u)$ is defined by $g'(u) = \displaystyle{-\frac{1}{2f'(u)}}$.
\end{thm}

\noindent {\it Proof:} It follows from \eqref{E:Eq-11} that the Gauss curvature $K = const \neq 0$  if and
only if the function $f(u)$ satisfies the following differential
equation
$$f''(u) + K f(u) = 0.$$
The general solution of the above equation is given by \eqref{E:Eq-Th}, 
where $\alpha$ and $\beta$ are constants.
The function $g(u)$ is
determined by $g'(u) = \displaystyle{-\frac{1}{2f'(u)}}$.

\qed

\section{Meridian surfaces  of parabolic type with constant mean curvature}

Surfaces with constant mean curvature in arbitrary spacetime
are important objects for the special role they play in the theory of
general relativity. The study of constant mean curvature surfaces
(CMC surfaces) involves not only geometric methods but also PDE
and complex analysis, that is why the theory of CMC surfaces is of
great interest not only for mathematicians but also for physicists
and engineers. Surfaces with constant mean curvature in
Minkowski space have been studied intensively in  the last years. See for example \cite{Bran}, \cite{Chav-Can}, 
\cite{Liu-Liu-1}, \cite{Lop-2}, \cite{Sa}. 

Let $\mathcal{M}'''_m$  be a meridian surface of parabolic type.
Equality \eqref{E:Eq-12} implies that the mean curvature of $\mathcal{M}'''_m$  is given by
\begin{equation} \label{E:Eq13}
|| H || = \sqrt{\frac{\varepsilon (\kappa^2 f'^2 - (f f'' + f'^2)^2)}{4f^2 f'^2}}.
\end{equation}

The following theorem gives the classification of the  meridian surfaces of parabolic type with constant mean curvature.

\begin{thm}  \label{T:mean curvature}
Let $\mathcal{M}'''_m$  be a meridian surface of parabolic type from the general class. Then
$\mathcal{M}'''_m$  has constant mean curvature $|| H || = a = const$, $a \neq 0$ if and only if
$\kappa = const = b, \; b \neq 0$,
and the meridian $m$ is determined by $f' = y(f)$
where
\begin{equation} \notag
\begin{array}{ll}
\vspace{2mm}
y(t) = \displaystyle{\frac{1}{t}\left(C \pm \frac{t}{2} \sqrt{b^2 - 4 a^2 t^2} \pm \frac{b^2}{4a}  \arcsin \frac{2at}{b}\right)}, \; C = const, & \quad \textrm{if} \quad \langle H, H \rangle >0, \\
\vspace{2mm}
y(t) = \displaystyle{\frac{1}{t}\left(C \pm \frac{t}{2} \sqrt{b^2 + 4 a^2 t^2} \pm \frac{b^2}{4a}  \ln |2 at + \sqrt{b^2 + 4 a^2 t^2} | \right)}, \; C = const, & \quad \textrm{if} \quad \langle H, H \rangle <0, \\
\end{array}
\end{equation}
$g(u)$ is defined by $g'(u) = \displaystyle{-\frac{1}{2f'(u)}}$.
\end{thm}

\noindent {\it Proof:}
Using \eqref{E:Eq13} we obtain that $||H|| = a$ if
and only if
$$\kappa^2(v) = \frac{(f f'' + f'^2)^2 + \varepsilon 4 a^2 f^2 f'^2}{f'^2},$$
which implies
\begin{equation} \label{E:Eq15}
\begin{array}{l}
\vspace{2mm}
\kappa = const = b, \; b \neq 0;\\
\vspace{2mm} (f f'' + f'^2)^2 + \varepsilon 4 a^2 f^2 f'^2 = b^2 f'^2.
\end{array}
\end{equation}
If we set $f' = y(f)$ in the second equality of \eqref{E:Eq15}, we obtain
that the function $y = y(t)$ is a solution of the following
differential equation
\begin{equation} \label{E:Eq16}
t y y' + y^2= \pm y \sqrt{b^2 - \varepsilon 4 a^2 t^2}.
\end{equation}
In the case $\varepsilon = 1$ the general solution of equation \eqref{E:Eq16} is given by the formula
\begin{equation} \label{E:Eq17}
y(t) = \displaystyle{\frac{1}{t}\left(C \pm \frac{t}{2} \sqrt{b^2 - 4 a^2 t^2} \pm \frac{b^2}{4a}  \arcsin \frac{2at}{b}\right)}, \qquad C = const.
\end{equation}
In the case $\varepsilon = -1$ the general solution of \eqref{E:Eq16} is given by
\begin{equation} \label{E:Eq17-a}
y(t) = \displaystyle{\frac{1}{t}\left(C \pm \frac{t}{2} \sqrt{b^2 + 4 a^2 t^2} \pm \frac{b^2}{4a}  \ln |2 at + \sqrt{b^2 + 4 a^2 t^2} | \right)}, \qquad C = const.
\end{equation}
The function $f(u)$ is determined by $f' = y(f)$ and \eqref{E:Eq17} or \eqref{E:Eq17-a}, respectively.
The function $g(u)$ is defined by $g'(u) = \displaystyle{-\frac{1}{2f'(u)}}$.

\qed

\section{Meridian surfaces  of parabolic type with constant invariant $k$}

Let $\mathcal{M}'''_m$ be a  meridian surface  of parabolic type.
Then the invariant $k$ is given by the formula
\begin{equation} \label{E:Eq19}
k = - \frac{\kappa_m^2(u) \, \kappa^2(v)}{f^2(u)}.
\end{equation}

In the following theorem we describe the meridian surfaces of parabolic type with constant invariant $k$.

\begin{thm} \label{T:constant k}
Let $\mathcal{M}'''_m$ be a meridian surface of parabolic type from the general class.
Then $\mathcal{M}'''_m$  has constant invariant $k = const = - a^2, \; a\neq 0$ if and only if  $\kappa = const = b, \; b \neq 0$,
and the meridian $m$ is determined by $f' = y(f)$
where
\begin{equation} \notag
y(t) = c \pm \frac{at^2}{2b}, \qquad c = const,
\end{equation}
$g(u)$ is defined by  $g'(u) = \displaystyle{-\frac{1}{2f'(u)}}$.
\end{thm}

\noindent {\it Proof:}
Using that $\kappa_m(u) = \displaystyle{\frac{f''(u)}{f'(u)}}$, from  \eqref{E:Eq19} we obtain that
$k = const = - a^2, \; a \neq 0$ if and only if
$$\kappa^2(v) = \frac{a^2 f^2(u) f'^2(u)}{f''\,^2(u)}.$$
The last equality implies
\begin{equation}\label{E:Eq20}
\begin{array}{l}
\vspace{2mm}
\kappa = const = b, \; b \neq 0;\\
\vspace{2mm}
 b f''(u) = \pm a f(u) f'(u).
\end{array}
\end{equation}
Setting $f' = y(f)$ in the second equality of \eqref{E:Eq20}, we obtain that
the function $y = y(t)$ is a solution of the following differential equation
\begin{equation} \notag
b y y'  = \pm aty.
\end{equation}
The general solution of the above equation is given by
\begin{equation}  \label{E:Eq21}
y(t) = c \pm \frac{at^2}{2b}, \qquad c = const.
\end{equation}
The function $f(u)$ is determined by $f' = y(f)$ and \eqref{E:Eq21}.
The function $g(u)$ is defined by   $g'(u) = \displaystyle{-\frac{1}{2f'(u)}}$.

\qed

\section{Chen meridian surfaces of parabolic type}

In  \cite{GM5} we showed that a spacelike surface in $\mathbb{R}^4_1$  is a non-trivial  Chen surface if and only if the invariant function $\lambda$ is zero.
In the next theorem we give the classification of all Chen meridian surfaces of parabolic  type.

\begin{thm} \label{T:Chen}
Let $\mathcal{M}'''_m$ be a meridian surface of parabolic type from the general class. Then
$\mathcal{M}'''_m$ is  a Chen surface  if and only if
$\kappa = const = b, \; b \neq 0$,
and the meridian $m$ is determined by $f' = y(f)$
where
\begin{equation} \notag
y(t) =  \frac{1}{2 c\,t^{\pm1}}  \left(c^2 t^{\pm2} +b^2 \right), \qquad c = const \neq 0,
\end{equation}
$g(u)$ is defined by  $g'(u) = \displaystyle{-\frac{1}{2f'(u)}}$.

\end{thm}

\noindent {\it Proof:}
 It follows from \eqref{E:Eq7} that $\lambda = 0$ if
and only if
\begin{equation} \notag
\kappa^2(v) =  \frac{f'^4(u) - f^2(u) f''^2(u)}{f'^2(u)},
\end{equation}
which implies
\begin{equation}\notag
\begin{array}{l}
\vspace{2mm}
\kappa = const = b, \; b \neq 0;\\
\vspace{2mm}
f'^4(u) - f^2(u) f''^2(u) = b^2 f'^2(u).
\end{array}
\end{equation}
Hence,  the function $f(u)$ is a solution of the
following differential equation:
\begin{equation} \label{E:Eq22}
f f'' = \pm f' \sqrt{f'^2 - b^2}.
\end{equation}
Setting $f' = y(f)$ in equation \eqref{E:Eq22}, we obtain
that the function $y = y(t)$ is a solution of the equation:
\begin{equation} \label{E:Eq23}
 t y y' = \pm y \sqrt{y^2 - b^2}.
\end{equation}
Since $y \neq 0$
the last equation is equivalent to
\begin{equation} \label{E:Eq24}
\frac{y'}{\sqrt{y^2 - b^2}} = \pm \frac{1}{t}.
\end{equation}
Integrating both sides of \eqref{E:Eq24}, we get
\begin{equation} \notag
 y + \sqrt{y^2 - b^2} = c \, t^{\pm1}, \qquad c = const.
\end{equation}
 Hence, the general solution of differential equation \eqref{E:Eq23} is given by
\begin{equation} \notag
y(t) =  \frac{1}{2 c\,t^{\pm1}}  \left(c^2 t^{\pm2} +b^2 \right), \qquad c = const \neq 0.
\end{equation}

\qed

\section{Meridian surfaces with parallel normal bundle}

Surfaces with parallel normal bundle are characterized by the condition  $\beta_1 = \beta_2 =0$ \cite{GM-new}.
In this section we  describe the meridian surfaces of parabolic type with parallel normal bundle.

\begin{thm} \label{T:parallel}
Let $\mathcal{M}'''_m$ be a meridian surface of parabolic type from the general
class. Then $\mathcal{M}'''_m$ has parallel normal bundle if and only if
one of the following cases holds:

\hskip 6mm (a) the meridian $m$ is defined by
\begin{equation} \notag
\begin{array}{l}
\vspace{2mm}
f(u) = \pm (cu + d)^{\frac{1}{2}};\\
\vspace{2mm}
g(u) = \displaystyle{\mp \frac{2}{3 c^2} (cu + d)^{\frac{3}{2}} + a},
\end{array}
\end{equation}
where  $a$, $c$, and $d$ are constants;

\hskip 6mm (b)  $\kappa = const = b, \; b \neq 0$, and the meridian $m$ is
determined by $f' = y(f)$ where
\begin{equation} \notag
y(t) = \frac{c + at}{t}, \quad a = const \neq 0, \quad c = const,
\end{equation}
$g(u)$ is defined by  $g'(u) = \displaystyle{-\frac{1}{2f'(u)}}$.
\end{thm}

\noindent {\it Proof:}
Using formulas \eqref{E:Eq7} we get
that $\beta_1 = \beta_2 =0$ if and only if
\begin{equation} \label{E:Eq26}
\begin{array}{l}
\vspace{3mm}
 \displaystyle{ \kappa\, \frac{d}{du}\left( \frac{f f'' + f'^2}{f'} \right) -  \frac{d}{dv}(\kappa) \, \frac{f f'' + f'^2}{f f' \sqrt{\dot{\varphi} ^2 + \varphi^2}}  = 0};\\
\vspace{2mm}
\displaystyle{ \kappa\, \frac{d}{du}\left( \frac{f f'' + f'^2}{f'} \right) +  \frac{d}{dv}(\kappa) \, \frac{f f'' + f'^2}{f f' \sqrt{\dot{\varphi} ^2 + \varphi^2}}  = 0}.
\end{array}
\end{equation}
It follows from \eqref{E:Eq26} that there are  two possible  cases:

\vskip 1mm
Case (a): $f f'' + f'^2 = 0$.
The general solution of this differential equation is  $f(u) = \pm \sqrt{cu +d}$, $c = const$,  $d = const$.
Using that $g'(u) = \displaystyle{-\frac{1}{2f'(u)}}$, we get  $g' = \displaystyle{\mp \frac{\sqrt{cu +d}}{c}}$. Integrating both sides of the last equation
we obtain $g(u) = \displaystyle{\mp \frac{2}{3 c^2} (cu + d)^{\frac{3}{2}} + a}$, $a = const$.
Consequently, the meridian $m$ is defined as described in \emph{(a)}.

\vskip 1mm
Case (b):
$ \displaystyle{\frac{f f'' + f'^2}{f'}}  = a = const$, $a \neq 0$ and $\kappa = b = const$, $b \neq 0$.
In this case the meridian $m$ is determined by the following differential
equation:
\begin{equation} \label{E:Eq27}
f f'' + f'^2 = a f', \qquad a = const \neq 0.
\end{equation}

The solutions of differential equation \eqref{E:Eq27} can be found in the  following way. Setting $f' = y(f)$ in equation \eqref{E:Eq27}, we obtain
that the function $y = y(t)$ is a solution of the equation:
\begin{equation} \label{E:Eq28}
t y y' +y^2= a y.
\end{equation}
Since $y \neq 0$
the last equation is equivalent to the equation
\begin{equation} \notag
 y' + \frac{1}{t} y = \frac{a}{t},
\end{equation}
whose general solution is given by the formula
\begin{equation} \notag
y(t) = \frac{c + at}{t}, \quad a = const \neq 0, \quad c = const.
\end{equation}

\qed

 \vskip 3mm \textbf{Acknowledgements:}
Research partially supported by the National Science Fund,
Ministry of Education and Science of Bulgaria under contract
DFNI-I 02/14.

\end{document}